\documentclass{amsart}
\usepackage[margin=1in]{geometry}

\usepackage{color,amssymb, comment, mathrsfs,tikz,upgreek,contour, soul, colonequals, mathdots,fancyvrb, amsmath}
\usepackage{tikz-cd}
\usetikzlibrary{cd}
\usetikzlibrary{fit,shapes.geometric}
\usetikzlibrary{matrix}
\usetikzlibrary{arrows}
\usepackage{enumitem}
\usepackage{mathtools}
\usepackage[all]{xy}
\usepackage[mathlines,pagewise]{lineno}
\usepackage[colorlinks, linkcolor={blue}]{hyperref}
\usepackage{cleveref}
\allowdisplaybreaks

\newtheorem{theorem}{Theorem}[section]
\newtheorem*{maintheorem1}{Main Theorem 1}
\newtheorem*{maintheorem2}{Main Theorem 2}

\newtheorem*{corollary2}{Corollary}
\newtheorem{lemma}[theorem]{Lemma}

\newtheorem{corollary}[theorem]{Corollary}
\theoremstyle{definition}
\newtheorem*{ack}{Acknowledgements}

\usepackage[only,llbracket,rrbracket,llparenthesis,rrparenthesis]{stmaryrd} 
\usepackage{accsupp}

\theoremstyle{definition}

\theoremstyle{remark}
\newtheorem{remark}[theorem]{Remark}

\numberwithin{equation}{section}

\newcommand{\bbZ}{\mathbb{Z}}

\newcommand{\kk}{\Bbbk}

\newcommand{\xra}{\xrightarrow}
\newcommand{\id}{\operatorname{id}}

\newcommand{\Ext}{\operatorname{Ext}}
\newcommand{\bExt}{\overline{\operatorname{Ext}}{\vphantom E}}

\newcommand{\Hom}{\operatorname{Hom}}
\newcommand{\bHom}{\overline{\operatorname{Hom}}}
\newcommand{\gdim}{\operatorname{G-dim}}

\newcommand{\Tor}{\operatorname{Tor}}

\newcommand{\m}{\mathfrak{m}}

\newcommand{\depth}{\operatorname{depth}}

\makeatletter
\DeclareFontEncoding{LS2}{}{\@noaccents}
\makeatother
\DeclareFontSubstitution{LS2}{stix}{m}{n}

\DeclareSymbolFont{largesymbolsstix}{LS2}{stixex}{m}{n}

\DeclareMathDelimiter{\lbrbrak}{\mathopen}{largesymbolsstix}{"EE}{largesymbolsstix}{"14}
\DeclareMathDelimiter{\rbrbrak}{\mathclose}{largesymbolsstix}{"EF}{largesymbolsstix}{"15}

\crefname{diagram}{diagram}{diagrams}
\Crefname{diagram}{Diagram}{Diagrams}
\creflabelformat{diagram}{#2(#1)#3}

\renewcommand{\H}{\mathrm{H}}



\author[Ferraro]{Luigi Ferraro}
\address[Luigi Ferraro]{School of Mathematical and Statistical Sciences \\ University of Texas Rio Grande Valley \\ Edinburg, TX 78539, U.S.A}
\email[]{luigi.ferraro@utrgv.edu}

\begin{document}
\title{On a question of Takahashi concerning direct summands of syzygies}

\begin{abstract}
For a finitely generated module $M$ of finite Gorenstein dimension
over a commutative noetherian local ring $(R,\m,\kk)$, we prove
that the natural map
$\Ext_R(\kk,M)\longrightarrow\Ext_R(\kk,M/\m M)$
is zero unless $R$ is Gorenstein. More generally, the same
conclusion holds whenever $\Ext_R^i(M,R)=0$ for $i\gg0$.
This gives a Gorenstein analogue of a regularity criterion of Avramov and
Iyengar. The proof uses the right $\Ext_R(\kk,\kk)$-module
structure of bounded cohomology studied by the author in a previous work. As an application, we show
that $R$ is Gorenstein if a nonzero quotient of a finite direct
sum of syzygies of $\kk$ has finite Gorenstein dimension.
In particular, this answers affirmatively a question of
Takahashi concerning direct summands of syzygies of the
residue field.
\end{abstract}
\maketitle

\section{Introduction}

Throughout this paper, $(R,\m,\kk)$ is a commutative noetherian
local ring, and modules are finitely generated. We write
$\Omega_R^n\kk$ for the $n$th syzygy in a minimal free resolution
of $\kk$, with $\Omega_R^0\kk=\kk$.

The homological properties of a local ring can often be detected
on syzygies of its residue field. Dutta proved that $R$ is regular
if some syzygy of $\kk$ has a nonzero direct summand of finite
projective dimension, see \cite{dutta}. Avramov
generalized this result by proving that every nonzero quotient
of a finite direct sum of syzygies of $\kk$ has maximal projective
complexity and curvature \cite[Corollary 9]{extremal}. In
particular, finite projective dimension of such a quotient forces
$R$ to be regular, recovering Dutta's theorem. Finite complete
intersection dimension of such a quotient similarly forces $R$
to be a complete intersection.

In the same paper, Avramov proved that tightly embeddable modules
also have maximal projective complexity and curvature
\cite[Theorem 4]{extremal}. Here an $R$-module $L$ is tightly
embeddable if it admits an embedding into an $R$-module $N$ such
that
\[
\m L\subsetneq\m N\subseteq L\subseteq N.
\]
Thus finite projective dimension or finite complete intersection
dimension of a tightly embeddable module characterizes regularity
or the complete intersection property, respectively. The
Gorenstein version of this criterion is also known: if a tightly
embeddable module has finite Gorenstein dimension, then $R$ is
Gorenstein. This follows from the more general result for Burch
submodules proved by Ghosh and Saha
\cite[Theorem 5.1]{GhoshSahaBurch}. In unpublished joint work
\cite{AFI}, Avramov, Iyengar, and the author obtain a different
proof of this criterion using a decomposition of stable
cohomology.

For direct summands of syzygies, the corresponding Gorenstein
question was posed by Takahashi \cite[Question 1.0.2]{tak}: if
$\Omega_R^n\kk$ has a nonzero direct summand of finite Gorenstein
dimension for some $n\geq0$, must $R$ be Gorenstein? The question
has been answered affirmatively in the following cases:
\begin{enumerate}
\item $n\leq\depth R+2$, by
\cite[Theorem 6.5]{tak2006}; the case $n\leq2$ was proved
in \cite{tak}.
\item $R$ is Cohen--Macaulay of minimal multiplicity
\cite[Corollary 5.3]{GhoshSyzygies}.
\item $R$ is a quotient of a Cohen--Macaulay local ring of minimal
multiplicity by a regular sequence
\cite[Theorem 7.5]{GhoshPuthenpurakal}.
\item $R$ is Cohen--Macaulay and admits a nonzero totally reflexive
module of minimal multiplicity and finite complexity
\cite[Corollary 3.7]{GheibiTakahashi}. This also yields the case
$\m^3=0$.
\item $R$ admits a regular sequence
$x_1,\ldots,x_d\in\m\setminus\m^2$ with
$\m^3\subseteq(x_1,\ldots,x_d)$
\cite[Corollary 3.3]{GhoshShort}.
\end{enumerate}
In this paper, we give an affirmative answer in general. More
precisely, we prove a Gorenstein analogue of Avramov's result for
nonzero quotients of finite direct sums of syzygies. Our main tool
is a Gorenstein version of a result of Avramov and Iyengar, which
we establish using the right module structure of bounded
cohomology described in \cite{bimodule}.

For an $R$-module $M$, let
\[
\pi_M:M\longrightarrow M/\m M
\]
be the canonical surjection, and set
\[
\varepsilon_M^n
\colonequals
\Ext_R^n(\kk,\pi_M):
\Ext_R^n(\kk,M)
\longrightarrow
\Ext_R^n(\kk,M/\m M).
\]
Avramov and Iyengar proved that if $M$ has finite projective
dimension and $\varepsilon_M^n\neq0$ for some $n$, then $R$ is
regular \cite[Theorem 2.4 (iii)$\Rightarrow$(i)]{AIBass}. Our first main result gives
a Gorenstein counterpart of this criterion. We write
$\bExt_R^i(M,N)$ for bounded cohomology, recalled in Section~2.

\begin{maintheorem1}[\Cref{thm:Gor}]
Let $M$ be an $R$-module such that $\varepsilon_M^n\neq0$ for some
$n\geq0$. If
\[
\bExt_R^i(M,\kk)=0
\]
for some $i\in\bbZ$ with $i+n>\depth R$, then $R$ is Gorenstein.
\end{maintheorem1}

In particular, the conclusion holds when $M$ has finite
Gorenstein dimension. More generally, it holds whenever
$\Ext_R^i(M,R)=0$ for $i\gg0$, since this implies
$\bExt_R^i(M,\kk)=0$ for $i\gg0$, see
\Cref{lem:vanishing}.

The proof uses the right multiplicative structure on bounded
cohomology studied by the author in \cite{bimodule}. By
\cite[Corollary 3.5]{bimodule}, there is an isomorphism
\[
\bExt_R(\kk,\kk)
\cong
\Ext_R(\kk,R)\otimes_\kk\Tor^R(\kk,\kk)
\]
of $\Ext_R(\kk,\kk)$-bimodules. The right action on the tensor
product uses the Hopf algebra structure of $\Ext_R(\kk,\kk)$.
We use this description to prove that if
$0\neq\eta\in\Ext_R^n(\kk,\kk)$ and
$\bExt_R^i(\kk,\kk)\eta=0$ for some $i+n>\depth R$, then $R$
is Gorenstein, see \Cref{lem:Gor}. The right action allows the multiplication maps
in the proof of Main Theorem 1 to factor through $\bExt_R^i(M,\kk)$ and an application of \Cref{lem:Gor} yields the desired result.

Our next main result applies this criterion to quotients of
syzygies and gives the Gorenstein analogue of
\cite[Corollary 9]{extremal}.

\begin{maintheorem2}[Theorem~3.8]
Let $M$ be a nonzero quotient of a finite direct sum of syzygies
of $\kk$. The following conditions are equivalent:
\begin{enumerate}
\item $R$ is Gorenstein.
\item $\gdim_R M<\infty$.
\item $\Ext_R^i(M,R)=0$ for $i\gg0$.
\item $\bExt_R^i(M,\kk)=0$ for $i\gg0$.
\end{enumerate}
\end{maintheorem2}

The connection with Main Theorem~1 is that
$\varepsilon_M^n\neq0$ precisely when there is a map
$\Omega_R^n\kk\to M$ whose image is not contained in $\m M$, see \Cref{lem:epsilon}.
For a nonzero quotient of a finite direct sum of syzygies, one
of the component maps has this property by Nakayama's lemma.
Since every direct summand is a quotient, we obtain the
following affirmative answer to Takahashi's question.

\begin{corollary2}[Corollary~3.9]
If a nonzero direct summand of a syzygy of $\kk$ has finite
Gorenstein dimension, then $R$ is Gorenstein.
\end{corollary2}

The paper is organized as follows. Section~2 recalls bounded
cohomology, the Hopf algebra structure on $\Ext_R(\kk,\kk)$,
and the right actions used in the proofs. Section~3 establishes
the Gorenstein criterion and its applications to quotients and
direct summands of syzygies.

\begin{ack}
\Cref{lem:vanishing} appears in the unpublished manuscript \cite{AFI} that the author wrote in 2013 with L. Avramov and S. Iyengar. The author was partly supported by the Simons Foundation grant MPS-TSM-00007849. 
\end{ack}

\section{Background and setup}
Throughout, $(R,\m,\kk)$ is a commutative noetherian local ring, and modules are finitely generated. For a module $M$ we set
\[
\pi_M:M\rightarrow M/\m M,\quad \varepsilon_M^n\colonequals\Ext_R^n(\kk,\pi_M):\Ext_R^n(\kk,M)\rightarrow\Ext_R^n(\kk,M/\m M).
\]
We use homological grading for complexes, so differentials have degree $-1$.
When using cohomological notation we set $X^i=X_{-i}$ and
$\H^i(X)=\H_{-i}(X)$. We write $|\eta|=i$ if
$\eta\in\Ext_R^i(M,N)$, while $|t|=-n$ if
$t\in\Tor_n^R(M,N)$.

\subsection*{Bounded cohomology}
Let $P\rightarrow M$ and $Q\rightarrow N$ be free resolutions of finitely generated modules. In cohomological notation the Hom complex has components
\[
\Hom_R(P,Q)^n=\prod_{j\in\bbZ}\Hom_R(P_j,Q_{j-n}),
\]
with differential
$\partial(f)=\partial^Qf-(-1)^nf\partial^P$ for
$f\in\Hom_R(P,Q)^n$. We call $f$ bounded if $f_j=0$ for all sufficiently
large $j$. These maps form a subcomplex with components
\[
\bHom_R(P,Q)^n=\bigoplus_{j\in\bbZ}\Hom_R(P_j,Q_{j-n}).
\]
The bounded cohomology of $M$ and $N$ is
\[
\bExt_R^n(M,N)\colonequals\H^n(\bHom_R(P,Q)),\qquad n\in\bbZ.
\]
It is independent of the chosen resolutions; see
\cite[Section 1.3]{stableExt}. The inclusion of bounded maps into all maps
induces the natural map
$\bExt_R^n(M,N)\rightarrow\Ext_R^n(M,N)$.

Products in cohomology are induced by composition. Composing a bounded map on either side with a
homogeneous map gives a bounded map. Consequently, composition induces
associative products
\[
\begin{aligned}
\Ext_R^p(N,L)\times\bExt_R^q(M,N)
&\longrightarrow\bExt_R^{p+q}(M,L),\\
\bExt_R^p(N,L)\times\Ext_R^q(M,N)
&\longrightarrow\bExt_R^{p+q}(M,L),
\end{aligned}
\]
where $L$ is another $R$-module. In particular, $\bExt_R(M,N)$ is an
$\Ext_R(N,N)$-$\Ext_R(M,M)$-bimodule.

\subsection*{The Hopf algebra structure of Ext}
The Yoneda algebra $\Ext_R(\kk,\kk)$ is the
universal enveloping algebra of the homotopy Lie algebra $\pi(R)$; see
\cite[Theorem 10.2.1]{IFR}. Therefore $\Ext_R(\kk,\kk)$ is a connected graded Hopf algebra. We write
\[
\Delta:\Ext_R(\kk,\kk)\longrightarrow\Ext_R(\kk,\kk)\otimes_\kk\Ext_R(\kk,\kk),\qquad
\epsilon:\Ext_R(\kk,\kk)\longrightarrow\kk,\qquad S:\Ext_R(\kk,\kk)\longrightarrow\Ext_R(\kk,\kk)
\]
for the coproduct, counit, and antipode, respectively. If $\xi\in\pi(R)$, then
\[
\Delta(\xi)=\xi\otimes1+1\otimes\xi,\qquad
\epsilon(\xi)=0,\qquad S(\xi)=-\xi.
\]

For homogeneous $\eta\in\Ext_R(\kk,\kk)$ we use the notation
\[
\Delta(\eta)=\sum\eta_{(1)}\otimes\eta_{(2)},
\qquad |\eta_{(1)}|+|\eta_{(2)}|=|\eta|,
\]
with homogeneous terms in the sum.

\subsection*{Right actions on Tor and bounded cohomology}
We regard $\Tor_q^R(\kk,N)$ as having cohomological degree $-q$.
Let $F\rightarrow\kk$ be a free resolution. If
$\eta\in\Ext_R^n(\kk,\kk)$ is represented by a cocycle
$a\in\Hom_R(F,F)^n$, then $a\otimes\id_N$ induces a map
\[
\eta\cdot-:\Tor_q^R(\kk,N)\longrightarrow\Tor_{q-n}^R(\kk,N).
\]
This defines a left $\Ext_R(\kk,\kk)$-action by using the resolution of the first
argument $\kk$. We use the antipode to turn it into a right action:
\[
t\eta\colonequals(-1)^{|t||\eta|}S(\eta)t.
\]
This is the right action used in \cite{bimodule}.

The right $\Ext_R(\kk,\kk)$-action on $\bExt_R(\kk,N)$ is induced by
precomposition. More explicitly, if $G\rightarrow N$ is a free
resolution and $\beta\in\bHom_R(F,G)^i$ represents $b$, then
\[
b\eta=[\beta a]\in\bExt_R^{i+n}(\kk,N),
\]
where $a\in\Hom_R(F,F)^n$ represents $\eta$. The same construction
gives the usual right $\Ext_R(\kk,\kk)$-action on $\Ext_R(\kk,R)$.

By \cite[Corollary 3.5]{bimodule}, there is an isomorphism of graded
right $\Ext_R(\kk,\kk)$-modules
\[
\bExt_R(\kk,N)\cong
\Ext_R(\kk,R)\otimes_\kk\Tor^R(\kk,N),
\]
where the action on the tensor product is
\[
(b\otimes t)\eta=
\sum(-1)^{|t||\eta_{(1)}|}
(b\eta_{(1)})\otimes(t\eta_{(2)}).
\]
In this formula the action on the first factor is Yoneda composition,
and the action on the second factor is the one defined above.
\subsection*{Soft truncation below}
Let $X$ be an $R$-complex and let $n$ be an integer. Following
\cite[Definition 2.5.23]{LarsBook}, we write $X_{\supseteq n}$ for the
soft truncation below of $X$ at homological degree $n$. Its components are
\[
(X_{\supseteq n})_j=
\begin{cases}
X_j & \text{if }j>n,\\
\ker\partial_n^X & \text{if }j=n,\\
0 & \text{if }j<n.
\end{cases}
\]
The differential is inherited from $X$. 

\subsection*{Gorenstein dimension}
For a module $G$, set $G^*\colonequals\Hom_R(G,R)$. The module $G$ is
totally reflexive if the natural map $G\rightarrow G^{**}$ is an
isomorphism and
\[
\Ext_R^i(G,R)=0=\Ext_R^i(G^*,R)\qquad\text{for all }i>0.
\]
For a nonzero module $M$, its Gorenstein dimension $\gdim_RM$ is the
least integer $s\geq0$ for which there is an exact sequence
\[
0\longrightarrow G_s\longrightarrow\cdots\longrightarrow G_0
\longrightarrow M\longrightarrow0
\]
with each $G_j$ totally reflexive, and is infinite if no such sequence
exists. If $\gdim_RM=s<\infty$, then
$\Ext_R^i(M,R)=0$ for all $i>s$. Moreover, $R$ is Gorenstein if and
only if every finitely generated module has finite Gorenstein dimension.
\section{Main Results}
We start by collecting a few preliminary lemmas needed in the proofs of our main results.
\begin{lemma}\label{lem:duality}
Let $0\neq\eta\in\Ext_R^n(\kk,\kk)$. Then there is $t\in\Tor_n^R(\kk,\kk)$ such that $t\eta=1$.
\end{lemma}
\begin{proof}
Let $F$ be an acyclic closure of $\kk$. By minimality
\[
\Ext_R^n(\kk,\kk)=\Hom_R(F_n,\kk),\quad\Tor_n^R(\kk,\kk)=\kk\otimes_RF_n.
\]
Let $S$ be the antipode map of $\Ext_R(\kk,\kk)$. Let $a:F_n\rightarrow\kk$ be a representative of $S(\eta)$, which is nonzero since $S$ is bijective. Let $x\in F_n$ such that $c\colonequals a(x)\neq0$. Set $t=(-1)^nc^{-1}(1\otimes x)\in\Tor_n^R(\kk,\kk)$. Then
\begin{align*}
t\eta&=(-1)^{|t||\eta|}S(\eta)t\\
&=c^{-1}S(\eta)(1\otimes x)\\
&=c^{-1}(1\otimes a(x))\\
&=1.\qedhere
\end{align*}
\end{proof}

\begin{lemma}\label{lem:Gor}
Let $0\neq\eta\in\Ext_R^n(\kk,\kk)$. If $\bExt_R^i(\kk,\kk)\eta=0$ for some $i$ with $i+n>\depth R$, then $R$ is Gorenstein.
\end{lemma}
\begin{proof}
We recall that
\begin{equation}\label{eq:bExt}
\bExt_R(\kk,\kk)\cong\Ext_R(\kk,R)\otimes_\kk\Tor^R(\kk,\kk)
\end{equation}
as $\Ext_R(\kk,\kk)$-bimodules by \cite[Corollary 3.5]{bimodule}.

By \Cref{lem:duality}, there is $t\in\Tor_n^R(\kk,\kk)$ such that $t\eta=1$. Consider the following maps
\[
\Ext_R^{i+n}(\kk,R)\rightarrow\bExt_R^{i}(\kk,\kk)\rightarrow\bExt_R^{i+n}(\kk,\kk)\rightarrow\Ext_R^{i+n}(\kk,R),
\]
where the map on the left is $b\mapsto b\otimes t$ via the isomorphism \eqref{eq:bExt}, the map in the middle is right multiplication by $\eta$, and the map on the right is the projection on $\Ext_R^{i+n}(\kk,R)\otimes_\kk\Tor^R_0(\kk,\kk)$ via \eqref{eq:bExt}. We prove that the composition of these three maps is the identity. Since the map in the middle is zero by hypothesis, it will follow that $\Ext_R^{i+n}(\kk,R)=0$, and therefore $R$ is Gorenstein by \cite[Proposition 2.3]{foxbyMu}.

Since $\Ext_R(\kk,\kk)$ is a connected graded Hopf algebra with coproduct $\Delta$ and counit $\epsilon$, it follows that $(\epsilon\otimes id)\Delta=id$, therefore
\[
\Delta(\eta)=1\otimes\eta+\sum_{|\eta_{(1)}|>0}\eta_{(1)}\otimes\eta_{(2)}.
\]
It follows that
\begin{align*}
(b\otimes t)\eta&=b\otimes(t\eta)+\sum_{|\eta_{(1)}|>0}(-1)^{|t||\eta_{(1)}|}(b\eta_{(1)})\otimes(t\eta_{(2)})\\
&=b\otimes1+\sum_{|\eta_{(1)}|>0}(-1)^{|t||\eta_{(1)}|}(b\eta_{(1)})\otimes(t\eta_{(2)}),
\end{align*}
where the second equality follows from the choice of $t$. Since $|(b\eta_{(1)})\otimes(t\eta_{(2)})|=i+n$ and $|b\eta_{(1)}|=i+n+|\eta_{(1)}|$, it follows that $t\eta_{(2)}\in\Tor^R_{|\eta_{(1)}|}(\kk,\kk)$. Therefore, projecting the expression in the display above onto $\Ext_R^{i+n}(\kk,R)\otimes_\kk\Tor^R_0(\kk,\kk)$ yields $b$, proving the claim.
\end{proof}
\begin{remark}
\Cref{lem:Gor} can also be proved under the hypothesis $\eta\bExt_R^i(\kk,\kk)=0$ for some $i$ with $i+n>\depth R$, with a simpler proof. But, for the purpose of this work we need the right $\Ext_R(\kk,\kk)$-module action on bounded cohomology.
\end{remark}
The next lemma in particular holds if $M$ has finite Gorenstein dimension.
\begin{lemma}\label{lem:vanishing}
If $\Ext_R^i(M,R)=0$ for $i>s$ for some $s\geq0$, then $\bExt_R^i(M,\kk)=0$ for $i>s$.
\end{lemma}

\begin{proof}
Let $P$ and $F$ be minimal free resolutions of $M$ and $\kk$ respectively. Denote $\Hom_R(P,R)$ by $P^*$. By \cite[Lemma 1.3.3]{stableExt} there is an isomorphism
\[
P^*\otimes_RF\cong\bHom_R(P,F),
\]
as DG $\Hom_R(F,F)$-$\Hom_R(P,P)$-bimodules.

Let $C\colonequals(P^*)_{\supseteq-s}$ be the soft truncation below of $P^*$ at homological degree $-s$. By \cite[Proposition 4.2.4]{LarsBook} and the hypothesis, the complexes $C$ and $P^*$ are quasi-isomorphic, and since $F$ is semi-free, there is a quasi-isomorphism
\[
C\otimes_RF\simeq P^*\otimes_RF.
\]
If $i>s$ then
\[
(C\otimes_RF)^i=\bigoplus_{p+q=i}C^p\otimes_RF^q=0,
\]
where the last equality follows because $F^q=F_{-q}=0$ for
$q>0$, while for $q\leq 0$ we have $p=i-q\geq i>s$, and hence
$C^p=C_{-p}=0$ by the definition of $C$. Consequently,
\[
\bExt_R^i(M,k)
\cong H^i(P^*\otimes_R F)
\cong H^i(C\otimes_R F)
=0
\qquad\text{for all }i>s.\qedhere
\]
\end{proof}
The result below can be interpreted as a Gorenstein version of \cite[Theorem 2.4 (iii)$\Rightarrow$(i)]{AIBass} since the finiteness of the Gorenstein dimension of $M$ implies the eventual vanishing of $\Ext_R^i(M,R)$, which by \Cref{lem:vanishing} implies the eventual vanishing of $\bExt_R^i(M,\kk)$.

\begin{theorem}\label{thm:Gor}
If $\varepsilon_M^n\neq0$ for some $n\geq0$, and $\bExt_R^{i}(M,\kk)=0$ for some $i$ with $i+n>\depth R$, then $R$ is Gorenstein.
\end{theorem}
\begin{proof}
Let $\alpha\in\Ext_R^n(\kk,M)$ such that $\varepsilon_M^n(\alpha)\neq0$. Fix a $\kk$-basis of $M/\m M$ and let $\lambda_j:M/\m M\rightarrow\kk$ for $j=1,\ldots,r$ be the coordinate functionals with respect to this basis. Then there is an isomorphism of $\kk$-vector spaces
\[
\begin{aligned}
\Ext_R^n(\kk,M/\m M)
&\xrightarrow{\ \cong\ }
\bigoplus_{j=1}^{r}\Ext_R^n(\kk,\kk),\\
y &\longmapsto
\bigl(\Ext_R^n(\kk,\lambda_j)(y)\bigr)_{j=1}^{r}.
\end{aligned}
\]
Since $\varepsilon_M^n(\alpha)\neq0$, there is a $j$ such that $\Ext_R^n(\kk,\lambda_j)(\varepsilon_M^n(\alpha))\neq0$. Set $\theta\colonequals \lambda_j\pi_M\in\Ext_R^0(M,\kk)$, then functoriality gives
\begin{samepage}
\begin{align*}
\Ext_R^n(\kk,\theta)(\alpha)&=\Ext_R^n(\kk,\lambda_j)(\Ext_R^n(\kk,\pi_M)(\alpha))\\
&=\Ext_R^n(\kk,\lambda_j)(\varepsilon_M^n(\alpha))\neq0.
\end{align*}
\end{samepage}
Set $\eta\colonequals\Ext_R^n(\kk,\theta)(\alpha)=\theta\alpha$, which is a nonzero element of $\Ext_R^n(\kk,\kk)$. Right multiplication by $\eta$ factors as
\begin{equation}\label{eq:factorization}
\bExt_R^i(\kk,\kk)\xra{\cdot\theta}\bExt_R^i(M,\kk)\xra{\cdot\alpha}\bExt_R^{i+n}(\kk,\kk),
\end{equation}
therefore the hypothesis implies $\bExt_R^i(\kk,\kk)\eta=0$, and \Cref{lem:Gor} implies that $R$ is Gorenstein.
\end{proof}
\begin{remark}
If one wishes to use the left $\Ext_R(\kk,\kk)$-action on $\bExt_R(\kk,\kk)$, then \eqref{eq:factorization} needs to be replaced with 
\[
\bExt_R^i(\kk,\kk)\xra{\alpha\cdot}\bExt_R^{i+n}(\kk,M)\xra{\theta\cdot}\bExt_R^{i+n}(\kk,\kk),
\]
but the eventual vanishing of the middle term is not guaranteed even when $M$ has finite Gorenstein dimension. Indeed, if $R$ is a ring that is not Gorenstein, and $M\colonequals R$, then
\[
\bExt^i_R(\kk,M)\cong\Ext_R^i(\kk,M)\neq0,\quad i\geq\depth R.
\]
On the other hand, the eventual vanishing of $\bExt^i_R(M,\kk)$ is guaranteed for modules of finite Gorenstein dimension by \Cref{lem:vanishing}. Since our main purpose is to answer Takahashi's question, this justifies the need to work with the right $\Ext_R(\kk,\kk)$-module structure of $\bExt_R(\kk,\kk)$.
\end{remark}

Before we answer Takahashi's question we prove one more lemma.

\begin{lemma}\label{lem:epsilon}
Let $n$ be a nonnegative integer. The following conditions are equivalent:
\begin{enumerate}
\item $\varepsilon_M^n\neq0$.
\item There is an $R$-linear map $\beta:\Omega_R^n\kk\rightarrow M$ such that $\beta(\Omega_R^n\kk)\not\subseteq\m M$.
\end{enumerate}
\end{lemma}
\begin{proof}
Let $F$ be a minimal free resolution of $\kk$. For an $R$-module $X$, we compute $\Ext_R(\kk,X)$ as the cohomology of $\Hom_R(F,X)$. We use $[-]$ to denote cohomology classes of maps in $\Hom_R(F,X)$. We note that since $F$ is minimal, the complex $\Hom_R(F,M/\m M)$ has zero differential. A cocycle of $\Hom_R(F,M)$ of cohomological degree $n$ is a map $p:F_n\rightarrow M$ with $p\partial^F_{n+1}=0$. Since $\Hom_R(F,M/\m M)$ has zero differential, it follows that
\[
\varepsilon_M^n([p])=[\pi_Mp]\neq0\Longleftrightarrow\pi_M p\neq0.
\]
We will also be denoting by $\rho_n$ the surjection $F_n\rightarrow\Omega_R^n\kk$.

$(1)\Rightarrow(2)$. Let $[p]\in\Ext_R^n(\kk,M)$ with $\varepsilon_M^n([p])\neq0$, represented by a map $p:F_n\rightarrow M$ with $p\partial^F_{n+1}=0$. Since $p\partial^F_{n+1}=0$, there is a map $\beta:\Omega_R^n\kk\rightarrow M$ such that $p=\beta\rho_n$. Since $\pi_Mp\neq0$, it follows that $\pi_M\beta\neq0$, and therefore $\mathrm{Im}\;\beta\not\subseteq \m M$.

$(2)\Rightarrow(1)$. Set $p\colonequals\beta\rho_n$. Since $\rho_n\partial_{n+1}^F=0$, the map $p$ is a cocycle of $\Hom_R(F,M)$ of cohomological degree $n$, and it therefore defines a class $[p]\in\Ext_R^n(\kk,M)$. By hypothesis $\pi_M\beta\neq0$, and the surjectivity of $\rho_n$ gives $\pi_Mp=\pi_M\beta\rho_n\neq0$. Therefore $\varepsilon_M^n([p])\neq0$.
\end{proof}

The following theorem can be interpreted as a Gorenstein version of \cite[Corollary 9]{extremal}.

\begin{theorem}
Let $M$ be a nonzero quotient of a finite direct sum of syzygies of $\kk$. The following are equivalent:
\begin{enumerate}
\item $R$ is Gorenstein;
\item $\gdim_RM<\infty$;
\item $\Ext_R^i(M,R)=0$ for $i\gg0$;
\item $\bExt_R^i(M,\kk)=0$ for $i\gg0$.
\end{enumerate}
\end{theorem}
\begin{proof}
The implications $(1)\Rightarrow(2)\Rightarrow(3)$ are standard results that hold for any module. The implication $(3)\Rightarrow(4)$ follows from \Cref{lem:vanishing} and also holds for any module. We prove $(4)\Rightarrow(1)$.

Let $\rho:\bigoplus_{a=1}^u\Omega_R^{n_a}\kk\rightarrow M$ be the given surjection. By Nakayama's Lemma, there is at least one component of $\rho$ with image not contained in $\m M$, denote it by $\beta:\Omega_R^{n_a}\kk\rightarrow M$. By \Cref{lem:epsilon} $\varepsilon_M^{n_a}\neq0$, and \Cref{thm:Gor} gives that $R$ is Gorenstein.
\end{proof}

As an immediate consequence we give an affirmative answer to a question of Takahashi \cite[Question 1.0.2]{tak}.

\begin{corollary}
If a nonzero direct summand of a syzygy of $\kk$ has finite Gorenstein dimension, then $R$ is Gorenstein.
\end{corollary}

\bibliographystyle{amsplain}
\bibliography{biblio}
\end{document}